\documentclass[12pt]{article}

\usepackage{amssymb}
\usepackage{amsmath}
\usepackage{amsthm}
\usepackage{xypic}

\renewcommand{\a}[1]{\ar@{-}[#1]}

\newtheorem{thm}{Theorem}[section]
\newtheorem{lem}[thm]{Lemma}

\newtheorem{cor}[thm]{Corollary}

\newtheorem*{thm*}{Theorem}

\theoremstyle{definition}

\theoremstyle{remark}

\newcommand{\system}[1]{\mbox{\fontfamily{cmss}\fontshape{n}\fontseries{m}%
    \selectfont#1}}
\newcommand{\RCA}{\system{RCA}\ensuremath{_0}}
\newcommand{\WKL}{\system{WKL}\ensuremath{_0}}
\newcommand{\ACA}{\system{ACA}\ensuremath{_0}}
\newcommand{\RTinffour}{\system{RT}\ensuremath{(4)}}
\newcommand{\RTfourk}{\system{RT}\ensuremath{^4_k}}
\newcommand{\RTnk}{\system{RT}\ensuremath{^n_k}}

\title{Reverse mathematics and infinite traceable graphs}

\author{Peter Cholak \\ \small{Department of Mathematics, University of
Notre Dame} \\ David Galvin \\ \small{Department of Mathematics,
University of Notre Dame} \\ Reed Solomon \\
\small{Department of Mathematics, University of Connecticut}}

\date{December 17, 2010}

\begin{document}

\maketitle

\footnotetext{Peter Cholak is partially supported by NSF grants
 DMS-0800198 and DMS-0652669 and David Galvin is supported in part by
 National Security Agency grant H98230-10-1-0364.}

\section{Introduction}

This paper falls within the general program of investigating the proof
theoretic strength (in terms of reverse mathematics) of combinatorial
principals which follow from versions of Ramsey's theorem.  We examine
two statements in graph theory and one statement in lattice theory
proved by Galvin, Rival and Sands \cite{GRS:82} using Ramsey's theorem
for 4-tuples.  Our main results are that the statements concerning
graph theory are equivalent to Ramsey's theorem for 4-tuples over
$\RCA$ while the statement concerning lattices is provable in $\RCA$.
We give the basic definitions for graph theory and lattice theory
below, but assume the reader is familiar with the general program of
reverse mathematics.  The definitions in this section are all given
within $\RCA$.

If $X \subseteq \mathbb{N}$ and $n \in \mathbb{N}$, then $[X]^n$
denotes the set of all $n$-element subsets of $X$.  A
$k$-\textit{coloring} of $[X]^n$ is a function $c:[X]^n \rightarrow
k$.  Ramsey's theorem for $n$-tuples and $k$ colors (denoted $\RTnk$)
is the statement that for all $k$-colorings of $[\mathbb{N}]^n$, there
is an infinite set $X$ such that $[X]^n$ is monochromatic.  Such a set
$X$ is called a \textit{homogeneous set} for the coloring.  We let
$\RTinffour$ denote the statement $\forall k \RTfourk$.  In terms of
reverse mathematics, $\RTnk$ is equivalent to $\ACA$ over $\RCA$ for
all $n \geq 3$ and $k \geq 2$ and $\RTinffour$ is equivalent to $\ACA$
over $\RCA$.  (See Section III.7 of Simpson \cite{sim:book}.)

Before giving the Galvin, Rival and Sands results, we introduce some
basic terminology from graph theory in $\RCA$.  A \textit{graph} $G$
is a pair $(V_G,E_G) = (V,E)$ such that $V$ (the \textit{vertex set})
is a subset of $\mathbb{N}$ and $E$ (the \textit{edge relation}) is a
symmetric irreflexive binary relation on $V$.  (Thus our graphs are
undirected and have no edges from a vertex to itself.)  If $E(x,y)$
holds, then we say there is an \textit{edge} between $x$ and $y$.
When specifying the edge relation on a graph, we assume that whenever
we say $E(x,y)$ holds we implicitly declare that $E(y,x)$ holds as
well.  (That is, we abuse notation by regarding $E(x,y)$ as shorthand
for $E(x,y) \wedge E(y,x)$.)  When we deal with more than one graph,
we denote the vertex set and edge relation of $G$ by $V_G$ and $E_G$.

A $n$--\textit{path} in a graph $G$ is a sequence of distinct vertices
$v_0, v_1, \ldots, v_{n-1}$ such that $E(v_i,v_{i+1})$ holds for all
$i \leq n-2$.  (A \textit{finite path} is an $n$--path for some $n \in
\mathbb{N}$.)  Similarly, an \textit{infinite path} is a sequence of
distinct vertices $v_0, v_1, \ldots$ (formally, specified by a
function $f:\mathbb{N} \rightarrow V$) such that $E(v_i,v_{i+1})$ for
all $i \in \mathbb{N}$.  If a path (finite or infinite) satisfies
$E(v_i,v_j)$ if and only if $|i-j| = 1$, then we say the path is
\textit{chordless}.  That is, a chordless path is a sequence of
vertices $v_0, v_1, \ldots$ (possibly finite) in which the only edges
are between vertices of the form $v_i$ and $v_{i+1}$.  (We use the
terminology of a chordless path from Galvin, Rival and Sands, but such
a path is also called an \textit{induced path} in the literature.)  An
infinite graph $G$ contains \textit{arbitrarily long chordless paths}
if for each $n \in \mathbb{N}$, $G$ contains a chordless $n$--path.
Similarly, we say $G$ contains an \textit{infinite chordless path} if
$G$ contains an infinite path which is chordless.

An infinite graph $G = (V,E)$ is \textit{traceable} if there is a
bijection $T:\mathbb{N} \rightarrow V$ (called a \textit{tracing
  function}) such that for all $i \in \mathbb{N}$, $E(T(i),T(i+1))$.
Thus, a traceable graph is one in which there is a path containing all
the vertices.  (A similar definition can be given when $G$ is finite.)

A graph $G = (V,E)$ is \textit{bipartite} if there is a partition $V =
V_0 \cup V_1$ such that for each edge $E(x,y)$ there is an $i \in \{
0,1 \}$ such that $x \in V_i$ and $y \in V_{1-i}$.  We use three
specific bipartite graphs in this paper.  The first graph is $K_{2,2}$
which has four vertices $a_0$, $a_1$, $b_0$ and $b_1$ with edges
between $a_i$ and $b_j$ for $i,j \leq 1$.

\bigskip

$
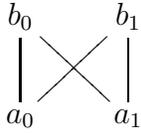
\xymatrix{ \a{d} b_0 \a{dr} & \a{d} b_1 \a{dl}  \\
  a_0 & a_1 }$

\bigskip
\noindent
The second graph is the complete countable bipartite graph
$K_{\omega,\omega}$.  Its vertices are $V = V_0 \cup V_1$ where $V_0 =
\{ a_n \mid n \in \mathbb{N} \}$ and $V_1 = \{ b_n \mid n \in
\mathbb{N} \}$ with edges between $a_n$ and $b_m$ for all $n,m \in
\mathbb{N}$.

\bigskip

$
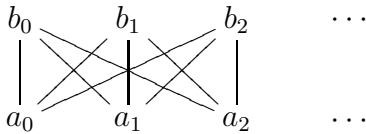
\xymatrix{ \a{d} b_0 \a{dr} \a{drr} & \a{d} b_1 \a{dl} \a{dr} & \a{d} b_2 \a{dl} \a{dll} & \cdots \\
  a_0 & a_1 & a_2 & \cdots }$

\bigskip
\noindent
Following the notation of \cite{GRS:82}, the third graph will be
denoted $A$.  Its vertices are $V = V_0 \cup V_1$ where $V_0 = \{ a_n
\mid n \in \mathbb{N} \}$ and $V_1 = \{ b_n \mid n \in \mathbb{N} \}$
with edges between $a_n$ and $b_m$ for all $n \leq m$.

\bigskip

$
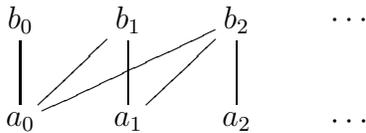
\xymatrix{ \a{d} b_0 & \a{d} b_1 \a{dl} & \a{d} b_2 \a{dl} \a{dll} & \cdots \\
  a_0 & a_1 & a_2 & \cdots }$

\bigskip

If $G$ and $H$ are graphs, then we say $G$ \textit{contains a subgraph
  isomorphic to} $H$ (or $G$ \textit{contains a copy of} $H$, or
\textit{there is an embedding of} $H$ \textit{into} $G$), if there is
an injective function $g: V_H \rightarrow V_G$ such that for all $x,y
\in V_H$, if there is an edge between $x$ and $y$ in $H$, then there
is an edge between $g(x)$ and $g(y)$ in $G$.  (Note that we allow
additional edges in $G$ between elements in the range of $g$.)  The
two graph theoretic results in \cite{GRS:82} are as follows.

\begin{thm}[Galvin, Rival and Sands \cite{GRS:82}]
  \label{thm:GRS1}
  Every infinite traceable graph either contains arbitrarily long
  finite chordless paths or contains a subgraph isomorphic to $A$.
\end{thm}

\begin{thm}[Galvin, Rival and Sands \cite{GRS:82}]
  \label{thm:GRS2}
  Every infinite traceable graph containing no chordless 4--path
  contains a subgraph isomorphic to $K_{\omega,\omega}$.
\end{thm}

As an application of Theorem \ref{thm:GRS1}, Galvin, Sands and Rival
prove the following lattice theoretic result.  (The lattice
terminology is defined in Section \ref{sec:lattices}.)

\begin{thm}[Galvin, Rival and Sands \cite{GRS:82}]
  \label{thm:GRS3}
  Every finitely generated infinite lattice of length 3 contains
  arbitrarily long finite fences.
\end{thm}

In Section \ref{sec:traceable}, we show that Theorems \ref{thm:GRS1}
and \ref{thm:GRS2} are equivalent to $\ACA$ over $\RCA$.  In Section
\ref{sec:lattices}, we show that Theorem \ref{thm:GRS3} is provable in
$\RCA$.  We follow Simpson \cite{sim:book} for the reverse mathematics
and we follow Soare \cite{soa:book} for computability theory.

\section{Traceable graphs}
\label{sec:traceable}

We begin this section with a computable combinatorics result which
will translate into a result in reverse mathematics.  If $G = (V,E)$
is a graph and $x \in V$, then we say $x$ has \textit{infinite degree}
if there are infinitely many $y$ such that $E(x,y)$.  Let $V^{\infty}$
denote the set of vertices with infinite degree in $G$.

\begin{thm}
  \label{thm:computable1}
  There is an infinite computable graph $G = (V,E)$ such that $G$ has
  a computable tracing function, $G$ has no chordless 4--paths and
  \[
  \forall X \bigg( \big( \exists e \, (W_e^X \text{is infinite} \,
  \wedge \, W_e^X \subseteq V^{\infty}) \big) \rightarrow 0' \leq_T X
  \bigg)
  \]
\end{thm}

\begin{proof}
  We build $G$ in stages using a dump construction to create a computable sequence of 
  nested subgraphs.  At stage $s$, our
  graph $G_s = (V_s,E_s)$ has $V_s = \{ 0, 1, \ldots, k_s \}$ for some
  $k_s \in \mathbb{N}$.  If $t > s$, then $k_t > k_s$ and for any $x,y
  \in V_s$, $E_t(x,y)$ holds if and only if $E_s(x,y)$ holds.  Thus,
  the vertex set $V = \cup_s V_s = \mathbb{N}$ is computable, the edge
  relation $E = \cup_s E_s$ is computable, and hence $G$ is
  computable.

  At stage $s$, the vertex set $V_s$ will be further subdivided into
  nonempty convex blocks $B_{0,s}, B_{1,s}, \ldots, B_{s,s}$.  (That
  is, if $n < p < m$ and $n,m \in B_{j,s}$, then $p \in B_{j,s}$.)
  Each block $B_{j,s}$ will have a coding vertex $c_{j,s}$, which will
  be the largest element of the block.  Thus, $B_{0,s} = \{ x \mid 0
  \leq x \leq c_{0,s} \}$ and for $0 < j \leq s$, $B_{j,s} = \{ x \mid
  c_{j-1,s} < x \leq c_{j,s} \}$.  For $x, y \leq k_s$, we say $x$ and
  $y$ are in the same $s$-block if $\exists j \leq s (x,y \in
  B_{j,s})$ and we say $x$ and $y$ are in different $s$-blocks
  otherwise.

  At stage $s+1$, we may collapse a final segment of these blocks by
  picking a value $0 \leq n \leq s$ and ``dumping" all the blocks
  currently after $B_{n,s}$ into $B_{n,s+1}$, i.e.~setting $\cup_{n
    \leq m \leq s} B_{m,s} \subseteq B_{n,s+1}$.  When we do this, we
  will redefine the coding vertices $c_{m,s+1}$ for $m \geq n$ to be
  new large numbers.  In the end, each coding vertex will have a limit
  $c_n = \lim_s c_{n,s}$ and each block will have a finite limiting
  block $B_n = \lim_s B_{n,s}$.  The limiting coding vertices will
  satisfy $c_0 < c_1 < c_2 < \cdots$.

  The only vertices with infinite degree will be the limiting $c_n$
  coding vertices.  Suppose $X$ can enumerate an infinite set of
  infinite degree vertices.  Then $X$ can compute an infinite set of
  infinite degree vertices in increasing order and hence $X$ can
  compute a function $f$ such that $f(n) \geq c_n$.  Therefore, it
  suffices to construct $G$ so that any function dominating the
  sequence $c_0, c_1, \ldots$ can compute $0'$.  The obvious way to do
  this is to make sure that $n \in K$ if and only if $n \in K_{c_n}$.
  The idea of the construction is to dump later blocks into $B_{n,s}$
  if $n$ enters $K_s$ and redefine $c_{n,s+1} \geq s$.

  Fix an enumeration $K_s$ of the halting problem $K$ such that
  exactly one number enters $K_s$ at each stage $s$.  Our construction
  proceeds in stages as follows.  At stage $0$, set $V_0 = \{ 0 \}$
  (and thus $k_0 = 0$), $E_0 = \emptyset$, $B_{0,0} = \{ 0 \}$ and
  $c_{0,0} = 0$.

  At stage $s+1$, check to see if the number $n$ entering $K$ at stage
  $s$ is large ($n > s$) or small ($n \leq s$).  If a number $n > s$
  enters $K_s$, then define $G_{s+1}$ as follows.  Let $V_{s+1} = V_s
  \cup \{ k_s+1 \}$.  (Recall that $k_s$ is the largest number in
  $V_s$.  Thus $k_{s+1} = k_s+1$.)  For each $j \leq s$, leave the
  blocks $B_{j,s+1} = B_{j,s}$ unchanged and leave the coding vertices
  $c_{j,s+1} = c_{j,s}$ unchanged.  Define a new block $B_{s+1,s+1} =
  \{ k_s+1 \}$ containing the newly added vertex and set its coding
  vertex $c_{s+1,s+1} = k_s+1$.  Add new edges between $c_{s+1,s+1}$
  and the other coding vertices $c_{j,s+1}$ for $j \leq s$ and end the
  stage.  (That is, let $E_{s+1}$ contain $E_s$ plus the edges
  $E_{s+1}(c_{j,s+1},c_{s+1,s+1})$ for each $j \leq s$.)

  If a number $n \leq s$ enters $K$ at stage $s$, then define
  $G_{s+1}$ as follows.  Let $u = (s+1) - n$.  Expand $V_s$ to
  $V_{s+1}$ by adding $u+1$ many new vertices $k_s+1, k_s+2, \ldots,
  k_s+u+1$.  (Thus $k_{s+1} = k_s+u+1$.)  For each $j < n$, leave the
  blocks $B_{j,s+1} = B_{j,s}$ and the coding vertices $c_{j,s+1} =
  c_{j,s}$ unchanged.  Dump the current later blocks and one
  additional element $k_s+1$ into $B_{n,s+1}$, and redefine the coding
  vertex $c_{n,s+1} = k_s+1$.  That is, set
  \[
  B_{n,s+1} = \bigcup_{n \leq m \leq s} B_{m,s} \cup \{ k_s+1 \}
  \text{ and } c_{n,s+1} = k_s+1.
  \]
  Use the remaining new elements $k_s+2, \ldots, k_s+u+1$ to define
  new single element blocks $B_{n+1,s+1}, \ldots, B_{s+1,s+1}$ with
  the single elements as the designated coding vertices.  That is, for
  each $1 \leq v \leq u$ set
  \[
  B_{n+v,s+1} = \{ k_s+v+1 \} \text{ and } c_{n+v,s+1} = k_s+v+1
  \]
  Expand $E_s$ to $E_{s+1}$ by adding new edges between each of the
  new coding vertices $c_{a,s+1}$ (for $n \leq a \leq s+1$) and all
  the other coding vertices $c_{b,s+1}$ (for $0 \leq b \leq s+1$ with
  $b \neq a$).  Also add edges $E_{s+1}(c_{n,s+1},x)$ for all $x \in
  B_{n,s+1}$ with $x \neq c_{n,s+1}$.  End the stage.

  This completes the construction of $G$.  As indicated above, $G$ is
  computable because there is an edge between $x$ and $y$ only if
  there is an edge between them at the first stage $s$ at which $x,y
  \in G_s$.  We check the remaining properties in a series of lemmas.

  \begin{lem}
    \label{lem:greatest}
    $\forall s \, \forall j \leq s \, \forall x \in B_{j,s} \, [x \leq
    c_{j,s} \, \wedge \, (x \neq c_{j,s} \rightarrow
    E_s(x,c_{j,s}))]$.
  \end{lem}

\begin{proof}
  This fact by induction on $s$.  For $s = 0$, it holds for $c_{0,0}$
  since $B_{0,0}$ is a singleton set.  For $s+1$, we split into cases
  depending on whether a small number enters $K_s$.  If not, then the
  property holds for $j < s+1$ by the induction hypothesis and the
  fact that the blocks and coding locations indexed by $j < s+1$ do
  not change.  It holds for $j = s+1$ since $B_{s+1,s+1}$ is a
  singleton set.

  If $n \leq s$ enters $K_s$, then the property holds for $j < n$ by
  the induction hypothesis and the fact that the blocks and coding
  vertices indexed by $j < n$ do not change.  It holds for $j = n$
  because $c_{n,s+1} = k_s+1$ is the largest element of $B_{n,s+1}$
  and we add edges at stage $s+1$ between this coding vertex and all
  the elements of $B_{n,s+1}$.  It holds for $n < j \leq s+1$ because
  each block $B_{j,s+1}$ is a singleton set.
\end{proof}

\begin{lem}
  \label{lem:codeconnection}
  $\forall s \, \forall i \neq j \leq s \, (E_s(c_{i,s},c_{j,s}))$.
\end{lem}

\begin{proof}
  This lemma follows by induction on $s$ since we add edges between
  any new or redefined coding vertices and all other coding vertices
  at each stage.
\end{proof}

\begin{lem}
  \label{lem:tracing}
  $\forall s \, \forall d < k_s \, (E_s(d,d+1))$.
\end{lem}

\begin{proof}
  This fact follows by induction on $s$.  For $s=0$, it is vacuously
  true since $k_0 = 0$.  For $s+1$, we split into cases depending on
  whether a small number enters $K_s$.  If not, then $k_{s+1} =
  k_s+1$.  By the induction hypothesis, we have $E_s(d,d+1)$, and
  hence $E_{s+1}(d,d+1)$, for all $d < k_s$.  It remains to show that
  $E_{s+1}(k_s, k_{s+1})$.  Since $k_s$ is the greatest element in
  $V_s$, it is also the greatest element in $B_{s,s}$.  By Lemma
  \ref{lem:greatest}, $c_{s,s} = k_s$ and hence by the construction
  $c_{s,s+1} = k_s$.  Since $c_{s+1,s+1} = k_{s+1} = k_s+1$ and since
  we add the edge $E_{s+1}(c_{s+1,s+1},c_{s,s+1})$, we have the edge
  $E_{s+1}(k_s,k_{s+1})$ as required.

  For the remaining case, suppose $n \leq s$ enters $K_s$ and hence
  $k_{s+1} = k_s+u+1$ where $u = (s+1)-n$.  By the induction
  hypothesis, $E_s(d,d+1)$, and hence $E_{s+1}(d,d+1)$, holds for all
  $d < k_s$.  It remains to show that $E_{s+1}(k_s+v,k_s+v+1)$ holds
  for all $v \leq u$.  First consider when $v=0$.  By construction,
  $k_s, k_{s+1} \in B_{n,s+1}$ and $c_{n,s+1} = k_{s+1}$.  Since we
  add edges from $c_{n,s+1}$ to each element of $B_{n,s+1}$, we have
  $E_{s+1}(k_s,k_{s+1})$ as required.  Finally, consider when $v > 0$.
  In this case, $c_{n+v-1,s+1} = k_s+v$ and $c_{n+v,s+1} = k_s+v+1$.
  By construction, we add an edge between these coding vertices at
  stage $s+1$ and hence have $E_{s+1}(k_s+v,k_s+v+1)$ as required.
\end{proof}

By Lemma \ref{lem:tracing}, $E(n,n+1)$ holds for all $n$.  Since $T(n)
= n$ is a bijection from $\mathbb{N}$ to $V = \mathbb{N}$, $T(n) = n$
is a computable tracing function for $G$.  We next show that $G$ has
no chordless 4--paths.  It suffices to show that each $G_s$ has no
chordless 4--paths.  We need two additional technical lemmas before
establishing this fact.  The first technical lemma says that whenever
we have an edge $E_s(x,y)$ with $x < y$, then either $x$ and $y$ are
in the same $s$-block or $x$ is a coding vertex $x = c_{j,s}$ for some
$j \leq s$.

\begin{lem}
  \label{lem:components}
  $\forall s \, \forall x < y \in G_s \, (E_s(x,y) \rightarrow [
  \exists j \leq s \, (x,y \in B_{j,s}) \vee \exists j \leq s \, (x =
  c_{j,s}))]$.
\end{lem}

\begin{proof}
  This lemma follows by induction on $s$.  For $s=0$, it is trivial
  since $|G_0| = 1$.  For $s+1$, assume that $x < y \in G_{s+1}$ and
  $E_{s+1}(x,y)$.  We need to show that either $x,y \in B_{j,s+1}$ for
  some $j \leq s+1$ or $x$ has the form $c_{j,s+1}$.  We split into
  cases depending on which (if any) of $x$ and $y$ are in $G_s$.  If
  $x \not \in G_s$, then by construction, $x = c_{j,s+1}$ for some $j$
  and we are done.

  If $x \in G_s$ and $y \not \in G_s$, then $y$ has the form
  $c_{j,s+1}$ for either $j = s+1$ (if no small number entered $K_s$)
  or for some $n \leq j \leq s+1$ (if $n \leq s$ entered $K_s$).  In
  the former case, by construction $E_{s+1}(x,y)$ implies that $x =
  c_{j,s+1}$ for some $j \leq s$ and we are done.  In the latter case,
  we split into cases depending on whether $y = c_{n,s+1}$ or $y =
  c_{j,s+1}$ for $j > n$.  If $y = c_{n,s+1}$, then $E_{s+1}(x,y)$
  implies that either $x \in B_{n,s+1}$ (and we are done since $y =
  c_{n,s+1} \in B_{n,s+1}$) or $x = c_{l,s+1}$ for some $l \neq j$
  (and we are done).  If $y = c_{j,s+1}$ for $j > n$, then
  $E_{s+1}(x,y)$ implies that $x = c_{l,s+1}$ for some $l \neq j$ (and
  we are done).

  Therefore, we are left with the case when $x,y \in G_s$.  Since $x,y
  \in G_s$ and $E_{s+1}(x,y)$, $E_s(x,y)$ must hold.  By the induction
  hypothesis, either $x, y \in B_{j,s}$ for some $j \leq s$ or $x$ has
  the form $c_{j,s}$.  If $x$ and $y$ are in the same $s$-block, then
  by construction they are in the same $(s+1)$-block.  (This block may
  or may not have the same index at stage $s+1$ depending on whether
  dumping occurred at stage $s+1$.)

  Therefore, we are left with the case when $x,y \in G_s$, $x$ and $y$
  are not in the same $s$-block and hence $x = c_{j,s}$ for some $j
  \leq s$.  If $c_{j,s} = c_{j,s+1}$, then $x = c_{j,s+1}$ and we are
  done.  Therefore, assume that $c_{j,s} \neq c_{j,s+1}$.  By the
  construction, this only occurs when a number $n \leq s$ enters $K_s$
  and $j \geq n$.  In this case, $x = c_{j,s}$ is dumped into
  $B_{n,s+1}$.  Since $x < y$ and $y \in G_s$, $y$ must also be dumped
  into $B_{n,s+1}$.  Hence, we have $x,y \in B_{n,s+1}$ and are done.
\end{proof}

Our second technical lemma says that whenever we have vertices $x < y$
which are connected but in different $s$-blocks, then $x$ is connected
to all of the elements in the $s$-block containing $y$.

\begin{lem}
  \label{lem:goup}
  The following statement holds for all stages $s$.  Let $x < y \in
  G_s$ with $E_s(x,y)$ and let $j \leq s$ be such that $y \in
  B_{j,s}$.  If $x \not \in B_{j,s}$, then $E_{s}(x,z)$ holds for all
  $z \in B_{j,s}$.
\end{lem}

\begin{proof}
  We prove this lemma by induction on $s$.  If $s=0$ then the
  statement holds trivially.  For $s+1$, first consider the case when
  no small number enters $K_s$. Let $j$ be such that $y \in B_{j,s+1}$
  and assume $x \not \in B_{j,s+1}$.  If $j \neq s+1$, then we are
  done because $E_s(x,z)$ (and hence $E_{s+1}(x,z)$) holds for all $z
  \in B_{j,s} = B_{j,s+1}$ by the induction hypothesis.  If $j = s+1$,
  then $y = c_{s+1,s+1}$ and $B_{j,s+1} = \{ y \}$, so again we are
  done.

  Second assume that $n \leq s$ enters $K_s$.  As above, let $j \leq
  s+1$ be such that $y \in B_{j,s+1}$ and assume $x \not \in
  B_{j,s+1}$.  If $j < n$, then as above (since $B_{j,s+1} = B_{j,s}$)
  we are done by the induction hypothesis.  If $j > n$, then as above
  (since $B_{j,s+1} = \{ y \}$) we are done trivially.  Therefore,
  assume that $j = n$.  In this case, $B_{n,s+1} = \cup_{n \leq l \leq
    s} B_{l,s} \cup \{ c_{n,s+1} \}$.  By Lemma \ref{lem:components},
  $x < y$ and $x \not \in B_{n,s+1}$ implies that $x = c_{i,s+1}$ for
  some $i < n$.  By construction, $c_{i,s+1} = c_{i,s}$, so $x =
  c_{i,s}$.  Therefore, for all $l$ such that $n \leq l \leq s$, we
  have $x < c_{l,s}$, $x \not \in B_{l,s}$ and $E_s(x,c_{l,s})$ holds.
  By the induction hypothesis, $E_s(x,z)$ (and hence $E_{s+1}(x,z)$)
  holds for all $z \in \cup_{n \leq l \leq s} B_{l,s}$.  Furthermore,
  by construction, $E_{s+1}(x,c_{n,s+1})$ holds completing this case.
\end{proof}

\begin{lem}
  $\forall s \, (G_s \text{ has no chordless 4--paths})$.
\end{lem}

\begin{proof}
  We proceed by induction on $s$.  For $s=0$, it follows trivially
  since $|G_0| = 1$.  For $s+1$, split into cases depending on whether
  a small number enters $K_s$.

  First, assume that no small number enters $K_s$ and assume for a
  contradiction that there is a chordless 4--path $x_0, x_1, x_2, x_3$
  in $G_{s+1}$.  By definition, we have $E_{s+1}(x_i,x_{i+1})$ for $i
  < 3$ and no other edges between these nodes (except those induced by
  symmetry).  By the induction hypothesis, at least one $x_i$ must lie
  outside $G_s$ and hence we have $x_i = k_s+1 = c_{s+1,s+1}$ for some
  $i \leq 3$.  We break into cases depending on which $x_i$ is equal
  to $c_{s+1,s+1}$.  Notice that if $x_0,x_1,x_2,x_3$ is a chordless
  4--path, then $x_3,x_2,x_1,x_0$ is also a chordless 4--path.
  Therefore, by symmetry, it suffices to show that we cannot have $x_0
  = c_{s+1,s+1}$ or $x_1 = c_{s+1,s+1}$.  (Recall that the elements in
  a path are required to be distinct.  We use this fact repeatedly
  without mention.)
  \begin{itemize}
  \item If $x_0 = c_{s+1,s+1}$, then by constuction $E_{s+1}(x_0,x_1)$
    implies that $x_1 = c_{l,s+1}$ for some $l \leq s$.  We break into
    subcases depending on the form of $x_2$.
    \begin{itemize}

    \item Suppose $x_2 < x_1$ and $x_2 \not \in B_{l,s+1}$.  By Lemma
      \ref{lem:components}, $x_2 = c_{j,s+1}$ for some $j < l$ and
      hence $E_{s+1}(x_0,x_2)$ for a contradiction.

    \item Suppose $x_2 \in B_{l,s+1}$ and consider the form of $x_3$.
      If $x_3 \in B_{l,s+1}$, then we have $E_{s+1}(x_1,x_3)$ for a
      contradiction.  If $x_3 < x_2$ and $x_3 \not \in B_{l,s+1}$,
      then $x_3 = c_{j,s+1}$ for some $j < l$ and we have
      $E_{s+1}(x_0,x_3)$ for a contradiction.  The remaining case,
      $x_3 > x_2$ and $x_3 \not \in B_{l,s+1}$ is not possible by
      Lemma \ref{lem:components} since $x_2 \in B_{l,s}$ but $x_2 \neq
      c_{l,s+1}$.

    \item Suppose $x_2 > x_1$ (so $x_2 \not \in B_{l,s+1}$) and
      consider the form of $x_3$.  If $x_3$ is in the same
      $(s+1)$-block as $x_2$, then since $E_{s+1}(x_1,x_2)$ holds, we
      have by Lemma \ref{lem:goup} that $E_{s+1}(x_1,x_3)$ holds for a
      contradiction.  If $x_3 < x_2$ and is not in the same
      $(s+1)$-block as $x_2$, then by Lemma \ref{lem:components}, $x_3
      = c_{i,s+1}$ for some $i$ and hence $E_{s+1}(x_0,x_3)$ holds for
      a contradiction.  If $x_3 > x_2$ and is not in the same
      $(s+1)$-block as $x_2$, then by Lemma \ref{lem:components}, $x_2
      = c_{i,s+1}$ for some $i$ and $E_{s+1}(x_0,x_2)$ holds for a
      contradiction.
    \end{itemize}

  \item If $x_1 = c_{s+1,s+1}$, then by the construction, $x_0 =
    c_{l,s+1}$ and $x_2 = c_{m,s+1}$ for some $l \neq m$.  By Lemma
    \ref{lem:codeconnection}, $E_{s+1}(x_0,x_2)$ holds for a
    contradiction.
  \end{itemize}

  Next assume that $n \leq s$ enters $K_s$ and $x_0,x_1,x_2,x_3$ is a
  chordless 4--path.  By the induction hypothesis, at least one of the
  $x_i$ is not in $G_s$ and hence must have the form $x_i = c_{j,s+1}$
  for some $n \leq j \leq s+1$.  If $x_i = c_{j,s+1}$ for $n < j \leq
  s+1$, then since $B_{j,s+1} = \{ c_{j,s+1} \}$, the same argument as
  in the previous case (when no small number enters $K_s$) suffices to
  derive a contradiction.  Therefore, we can assume without loss of
  generality that the chordless path is contained in $B_{0,s+1} \cup
  \cdots \cup B_{n,s+1}$ and that $x_i = c_{n,s+1}$ for some $i \leq
  3$.  By symmetry, it suffices to consider the cases when $x_0 =
  c_{n,s+1}$ and $x_1 = c_{n,s+1}$.  (Below, we frequently use without
  mention that none of the $x_i$ have the form $c_{l,s+1}$ for $l > n$
  and that if $x_i \in B_{n,s+1}$ and $x_j \not \in B_{n,s+1}$, then
  $x_j < x_i$.)

  \begin{itemize}

  \item Suppose $x_0 = c_{n,s+1}$ and consider the form of $x_1$.
    Since $x_1 < x_0$, either $x_1 \in B_{n,s+1}$ or $x_1 = c_{l,s+1}$
    for some $l < n$.  Consider these cases separately.

    \begin{itemize}

    \item Suppose $x_1 \in B_{n,s+1}$ and consider the form of $x_2$.
      If $x_2 \in B_{n,s+1}$, then $E_{s+1}(x_0,x_2)$ holds (since
      $x_0 = c_{n,s+1}$ is connected to all vertices in $B_{n,s+1}$)
      for a contradiction.  If $x_2 \not \in B_{n,s+1}$, then $x_2 <
      x_1$ and hence by Lemma \ref{lem:components}, $x_2 = c_{l,s+1}$
      for some $l < n$.  Thus $E_{s+1}(x_0,x_2)$ holds for a
      contradiction.

    \item Suppose $x_1 = c_{l,s+1}$ for some $l < n$ and consider the
      form of $x_2$.  There are three cases to consider.

      \begin{itemize}
      \item Assume $x_2 \in B_{l,s+1}$ and consider the form of $x_3$.
        If $x_3 \in B_{l,s+1}$, then by Lemma \ref{lem:greatest},
        $E_{s+1}(x_1,x_3)$ holds for a contradiction.  If $x_3 \not
        \in B_{l,s+1}$, then by Lemma \ref{lem:components} and the
        fact that $x_2 \in B_{l,s+1}$ but $x_2 \neq c_{l,s+1}$, we
        have $x_3 < x_2$ and hence $x_3 = c_{i,s+1}$ for some $i < l$.
        But then $E_{s+1}(x_0,x_3)$ holds for a contradiction.

      \item Assume $x_2 < x_1$ and $x_2 \not \in B_{l,s+1}$.  By Lemma
        \ref{lem:components}, $x_2 = c_{i,s+1}$ for some $i < l$ and
        hence $E_{s+1}(x_0,x_2)$ holds for a contradiction.

      \item Assume $x_2 > x_1$ (so $x_2 \not \in B_{l,s+1}$) and
        consider the form of $x_3$.  If $x_3$ is the same
        $(s+1)$-block as $x_2$, then since $E_{s+1}(x_1,x_2)$ holds,
        Lemma \ref{lem:goup} implies $E_{s+1}(x_1,x_3)$ holds for a
        contradiction.  Therefore, $x_3$ is not in the same
        $(s+1)$-block as $x_2$.  Therefore, by Lemma
        \ref{lem:components}, either $x_2$ or $x_3$ has the form
        $c_{i,s+1}$ for some $i$.  Hence either $E_{s+1}(x_0,x_2)$
        holds or $E_{s+1}(x_0,x_3)$ holds, giving a contradiction.
      \end{itemize}
    \end{itemize}

  \item Suppose $x_1 = c_{n,s+1}$.  By the construction, $x_1$ is
    connected only to the vertices in $B_{n,s+1}$ and the coding
    vertices $c_{l,s+1}$.  Since $E_{s+1}(x_0,x_1)$ and
    $E_{s+1}(x_1,x_2)$ hold, either $x_0 \in B_{n,s+1}$ or $x_0 =
    c_{l,s+1}$ for some $l < n$, and either $x_2 \in B_{n,s+1}$ or
    $x_2 = c_{i,s+1}$ for some $i < n$.  Consider each of the possible
    combinations separately.

    \begin{itemize}

    \item Suppose $x_0 = c_{l,s+1}$ and $x_2 = c_{i,s+1}$.  In this
      case, $E_{s+1}(x_0,x_2)$ holds for a contradiction.

    \item Suppose $x_0 = c_{l,s+1}$ and $x_2 \in B_{n,s+1}$.  Since
      $x_0 < x_1$, $x_0 \not \in B_{n,s+1}$, $E_{s+1}(x_0,x_1)$ holds
      and $x_1,x_2 \in B_{n,s+1}$, Lemma \ref{lem:goup} implies that
      $E_{s+1}(x_0,x_2)$ holds for a contradiction.

    \item Suppose $x_0 \in B_{n,s+1}$ and $x_2 = c_{i,s+1}$.  Since
      $x_2 < x_1$, $x_2 \not \in B_{n,s+1}$, $E_{s+1}(x_2,x_1)$ holds
      and $x_0,x_1 \in B_{n,s+1}$, Lemma \ref{lem:goup} implies that
      $E_{s+1}(x_2,x_0)$ holds for a contradiction.

    \item Suppose $x_0, x_2 \in B_{n,s+1}$.  Consider the form of
      $x_3$.  If $x_3 \in B_{n,s+1}$, then $E_{s+1}(x_1,x_3)$ holds
      for a contradiction. Therefore, $x_3 \not \in B_{n,s+1}$ and
      $x_3 < x_2$.  By Lemma \ref{lem:components}, $x_3 = c_{j,s+1}$
      for some $j < n$.  By construction $E_{s+1}(x_1,x_3)$ holds for
      a contradiction.

    \end{itemize}

  \end{itemize}
\end{proof}

We have now established that $G$ is a computable graph with a
computable tracing function and no chordless 4-paths.  It remains to
show that if $X$ can enumerate an infinite set of infinite degree
vertices, then $0' \leq_T X$.

\begin{lem}
  \label{lem:limits}
  $\forall k \, ( \lim_s c_{k,s} = c_k \text{ exists})$.
\end{lem}

\begin{proof}
  For any stage $s \geq k$, $c_{k,s+1} \neq c_{k,s}$ only if the block
  $B_{k,s}$ is dumped at stage $s+1$ into a block $B_{n,s+1}$ with $n
  \leq k$.  Since this happens only if a number $n \leq k$ enters
  $K_s$, we have that $c_{k,s}$ can change at most $k+1$ many times
  after it is first defined.
\end{proof}

From Lemma \ref{lem:limits} and the construction it is clear that for
all indices $k$ and all stages $s \geq k$, $c_{k,s} \leq c_{k,s+1}$.
Therefore, each $c_{k,s}$ is increasing in $s$ and stabilizes when it
reaches its limit.  It is also clear that $c_0 < c_1 < \cdots$ and
that $x \leq c_x$ for all $x$.  Finally, since for all stages $s$,
$B_{0,s} = \{ x \mid 0 \leq x \leq c_{0,s} \}$ and $B_{j,s} = \{ x
\mid c_{j-1,s} < x \leq c_{j,s} \}$ for $0 < j \leq s$, we have that
each block reaches a limit $B_j = \lim_s B_{j,s}$ and each vertex $x$
sits inside some limiting block.  (That is, for each vertex $x$, there
is a stage $s$ and a block $B_j$ such that $x \in B_{j,s} = B_j$.)

\begin{lem}
  \label{lem:infinitedegree}
  A vertex $x \in G$ has infinite degree if and only if $x = c_k$ for
  some $k$.
\end{lem}

\begin{proof}
  First, note that each vertex $c_k$ has infinite degree since
  $E(c_k,c_l)$ holds for all $l \neq k$.  (More formally, if $s$ and
  $t$ are stages such that $c_{k,s} = c_k$ and $c_{l,s} = c_l$, then
  by stage $u = \max \{ s,t \}$ we have added an edge $E_u(c_k,c_l)$.)

  Second, let $x$ be a vertex such that $x \neq c_k$ for all $k$.
  Suppose for a contradiction that $x$ has infinite degree.  Fix a
  stage $s$ and a block such that $x \in B_{j,s} = B_j$.  Since $x
  \neq c_j$ and both $B_{j,s}$ and $c_{j,s}$ have reached limits, it
  follows that for all stages $t \geq s$, $x \in B_{j,t}$ and $x \neq
  c_{j,t}$.  Since $x$ is assumed to have infinite degree, there must
  be a vertex $y > c_j$ and a stage $t > s$ such that $E_t(x,y)$
  holds.  By Lemma \ref{lem:greatest}, $y \not \in B_{j,t}$ and hence
  (since $x < y$ and $E_t(x,y)$) by Lemma \ref{lem:components}, $x =
  c_{l,t}$ for some $l$.  Since $x \in B_{j,t}$ we must have $x =
  c_{j,t}$ for a contradiction.
\end{proof}

In addition to having $x \leq c_x$, it is clear that $s \leq k_s$ for
all $s$.

\begin{lem}
  \label{lem:computes}
  $\forall x \, (x \in K \Leftrightarrow k \in K_{c_x})$.
\end{lem}

\begin{proof}
  Suppose $x$ enters $K$ at stage $s$.  If $s < x$, then $s < c_x$ and
  hence $x \in K_{c_x}$.  If $x \leq s$, then at stage $s+1$ of the
  construction, we dump later blocks into $B_{x,s+1}$ and set
  $c_{x,s+1} = k_s+1$.  Therefore, $c_x > s$ and hence $x \in
  K_{c_x}$.
\end{proof}

\begin{lem}
  \label{lem:mainprop}
  If $X$ can enumerate an infinite set of infinite degree vertices in
  $G$ then $0' \leq_T X$.
\end{lem}

\begin{proof}
  Define a function $f \leq_T X$ by setting $f(0) =$ the first
  infinite degree vertex enumerated by $X$ and $f(n+1) =$ the first
  infinite degree vertex $y$ enumerate by $X$ such that $y > f(n)$.
  By Lemma \ref{lem:infinitedegree}, $f$ has the property that $c_x
  \leq f(x)$ for all $x$ and hence $x \in K$ if and only if $x \in
  K_{f(x)}$.
\end{proof}

This completes the proof of Theorem \ref{thm:computable1}.

\end{proof}

Since the graph $G$ constructed in Theorem \ref{thm:computable1} is
traceable and has no chordless 4-paths, Theorems \ref{thm:GRS1} and
\ref{thm:GRS2} tell us that $G$ has subgraphs isomorphic to $A$ and to
$K_{\omega,\omega}$.  However, if $f$ is an embedding of either $A$ or
$K_{\omega,\omega}$ into $G$, then $f$ can enumerate an infinite set
of infinite degree vertices in $G$.  Therefore, by Theorem
\ref{thm:computable1}, $0' \leq_T f$ for any embedding of $A$ or
$K_{\omega,\omega}$ into $G$.  Thus we have the following corollary
concerning the lack of effectiveness of Theorems \ref{thm:GRS1} and
\ref{thm:GRS2}.

\begin{cor}
  There is a computable graph $G$ with a computable tracing function
  and no chordless 4-paths such that $0'$ is computable from any
  embedding of $A$ or $K_{\omega,\omega}$ into $G$.
\end{cor}

We next translate this result into the language of reverse
mathematics.

\begin{thm}[\RCA]
  \label{thm:traceable}
  The following are equivalent.
  \begin{enumerate}
  \item[(1).] Theorem \ref{thm:GRS1}.
  \item[(2).] Theorem \ref{thm:GRS2}.
  \item[(3).] \ACA.
  \end{enumerate}
\end{thm}

\begin{proof}
  The fact that (3) implies (1) and (2) follows immediately from the
  proofs given in \cite{GRS:82}.  (Our Theorem \ref{thm:GRS1} is Theorem 1 in 
  \cite{GRS:82} and our Theorem \ref{thm:GRS2} is Theorem 2 in \cite{GRS:82}.)  
  The proofs translate easily into proofs in
  $\RCA+\RTinffour$.  Since $\ACA \vdash \RTinffour$, this gives the
  desired implications.

  We prove (1) implies (3) and (2) implies (3) with essentially the
  construction given in the proof of Theorem \ref{thm:computable1}.
  For the remainder of this proof we work in $\RCA$.  Fix a 1-to-1
  function $f:\mathbb{N} \rightarrow \mathbb{N}$.  It suffices to
  construct a graph $G$ so that any embedding of $A$ or
  $K_{\omega,\omega}$ into $G$ yields a $\Delta^0_1$ definition of the
  range of $f$.

  We build a graph $G$ in stages as in the proof of Theorem
  \ref{thm:computable1}.  At stage $0$, set $V_0 = \{ 0 \}$ (so $k_0 =
  0$), $E_0 = \emptyset$, $B_{0,0} = \{ 0 \}$ and $c_{0,0} = 0$.  At
  stage $s+1$, let $n = f(s)$ and split into cases depending on
  whether $n > s$ or $n \leq s$.

  If $n > s$, then define $G_{s+1}$ as follows.  Let $V_{s+1} = V_s
  \cup \{ k_s+1 \}$ and $k_{s+1} = k_s+1$.  For each $j \leq s$, let
  $B_{j,s+1} = B_{j,s}$ and $c_{j,s+1} = c_{j,s}$.  Define
  $B_{s+1,s+1} = \{ k_{s+1} \}$ and $c_{s+1,s+1} = k_{s+1}$.  Expand
  $E_s$ to $E_{s+1}$ by adding edges between $c_{s+1,s+1}$ and each
  $c_{j,s+1}$ for $j \leq s$.

  If $n \leq s$, then let $u = (s+1)-n$ and define $G_{s+1}$ as
  follows.  Let $k_{s+1} = k_s +u+1$ and define $V_{s+1} = V_s \cup \{
  x \mid k_s < x \leq k_{s+1} \} = \{ x \mid x \leq k_{s+1} \}$.  For
  $j < n$, let $B_{j,s+1} = B_{j,s}$ and $c_{j,s+1} = c_{j,s}$.  Set
  \[
  B_{n,s+1} = \bigcup_{n \leq m \leq s} B_{m,s} \cup \{ k_s+1 \}
  \text{ and } c_{n,s+1} = k_s+1.
  \]
  For $1 \leq v \leq u$, set
  \[
  B_{n+v,s+1} = \{ k_s+v+1 \} \text{ and } c_{n+v,s+1} = k_s+v+1.
  \]
  Expand $E_s$ to $E_{s+1}$ by adding edges between each pair
  $c_{j,s+1}$ and $c_{i,s+1}$ with $i \neq j \leq s+1$.  (If $i,j < n$
  then these edges already exist in $E_s$.)

  Let $G = (V,E)$ where $V = \cup_s V_s = \mathbb{N}$ and $E = \cup_s
  E_s$.  Lemmas \ref{lem:greatest}, \ref{lem:codeconnection},
  \ref{lem:tracing}, \ref{lem:components} and \ref{lem:goup} were all
  proved by $\Sigma^0_1$ induction and hence are provable in $\RCA$.
  Therefore, $G$ is traceable and has no chordless 4-paths.

  We need an analog of Lemma \ref{lem:infinitedegree}.  Suppose $x \in
  G$ and $x$ is placed in $G$ at stage $s$.  By the construction, $x =
  c_{j,s}$ for some $j \leq s$.  If $\forall t > s \, (x = c_{j,t})$,
  then $x$ has infinite degree because we add an edge between $x$ and
  each new element added at stage $t$ for all $t > s$.  On the other
  hand, if $\exists t > s (x \neq c_{j,t})$, then by construction $x$
  is never equal to a coding vertex $c_{i,u}$ for any $u \geq t$.
  Since the only edges added at stages $u \geq t$ are between vertices
  of the form $c_{i,u}$ and $c_{j,u}$, $x$ is never connected by an
  edge to another vertex after stage $t$.  Therefore, $x$ has finite
  degree.

  We also need an analog of Lemma \ref{lem:computes}.  We claim that
  \[
  \forall k \, \forall s \geq k \, ( c_{k,s+1} \neq c_{k,s}
  \leftrightarrow f(s) \leq k).
  \]
  If $f(s) \leq k$, then the dumping in the definition of $G_{s+1}$
  causes $c_{k,s}$ to be redefined and hence $c_{k,s+1} \neq c_{k,s}$.
  On the other hand, if $c_{k,s+1} \neq c_{k,s}$ then dumping must
  have occurred because $f(s) = n \leq s$.  Furthermore, since
  $c_{j,s+1} = c_{j,s}$ for all $j < n$, we cannot have $k < n$.
  Therefore, $f(s) \leq k$.

  We are now ready to apply (1) or (2) and extract a definition of the
  range of $f$.  Since $G$ is traceable and has no chordless 4-paths,
  by (1) or (2), there is an embedding $g: A \rightarrow G$ or $g:
  K_{\omega,\omega} \rightarrow G$.  Recall that each of the vertices
  $a_k$ for $k \in \mathbb{N}$ in $A$ or $K_{\omega,\omega}$ have
  infinite degree.

  Define an auxiliary function $g': \mathbb{N} \rightarrow \mathbb{N}$
  by $g'(n) = \max \{ g(a_k) \mid k \leq n \}$.  (Note that the sets
  $\{ g(a_k) \mid k \leq n \}$ exist by bounded $\Sigma^0_1$
  comprehension.)  We claim that
  \[
  \forall k \, \forall t \geq g'(k) \, (c_{k,g'(k)} = c_{k,t}).
  \]
  Suppose for a contradiction that this property fails for some $k$
  and fix $t \geq g'(k)$ such that $c_{k,g'(k)} \neq c_{k,t}$.  Let $x
  = c_{k,g'(k)}$.  By our analog of Lemma \ref{lem:infinitedegree},
  $x$ has finite degree.  However, by the definition of $g'(k)$, $x =
  g(a_i)$ for some $i \leq k$.  Thus $a_i$ has infinite degree in $A$
  or $K_{\omega,\omega}$ but $g(a_i)$ has finite degree in $G$,
  contradicting the fact that $g$ is an embedding.

  By our analog of Lemma \ref{lem:computes}, this property implies
  that
  \[
  \forall t \geq g'(k) \, ( f(t) \not \leq k).
  \]
  Thus, $k$ is in the range of $f$ if and only if $\exists x \leq
  g'(k) \, (f(x) = k)$, completing the proof that (1) and (2) imply
  (3).
\end{proof}

\section{Finitely generated lattices}
\label{sec:lattices}

Our goal for this section is to show that Theorem \ref{thm:GRS3} is
provable in $\RCA$ and hence its proof does not require the use of
Theorem \ref{thm:GRS1}.  Before giving the formal lattice theoretic
definitions, we prove a finite Ramsey style result, Theorem
\ref{thm:kramsey} below, that is contained in \cite{GRS:82}.  The
proof of Theorem \ref{thm:kramsey} given in \cite{GRS:82} explicitly
uses Theorem \ref{thm:GRS1} (and hence this proof requires $\ACA$)
although the authors indicate that an alternate proof is available
using the Finite Ramsey Theorem.  Our proof of Theorem
\ref{thm:kramsey} formalizes this alternate approach in $\RCA$.

Recall that for $m,n,u,k \in \mathbb{N}$, the notation $[0,m]
\rightarrow (q)^n_k$ means that for any $k$-coloring of $[Y]^n$, where
$Y = [0,m]$, there is a set $X \subseteq Y$ such that $|X| = q$ and
the coloring is monochromatic on $[X]^n$.  The Finite Ramsey Theorem
is the statement
\[
\forall n,u,k \, \exists m \, ([0,m] \rightarrow (q)^n_k).
\]
The least $m$ satisfying $[0,m] \rightarrow (q)^n_k$ is called the
finite Ramsey number for $n$-tuples with $k$ many colors and a
homogeneous set of size $q$.  Since the Finite Ramsey Theorem is
provable in $PA^{-} + I\Sigma_1$ (see H\'ajek and Pudl\'ak
\cite{HP:book} Chapter II, Theorem 1.10) and $PA^{-} + I\Sigma_1$ is
the first order part of $\RCA$ (see Simpson \cite{sim:book} Corollary
IX.1.11), it follows that $\RCA$ proves the Finite Ramsey Theorem.
The finite style Ramsey result from \cite{GRS:82} is as follows.

\begin{thm}[$\RCA$]
  \label{thm:kramsey}
  For all $n \in \mathbb{N}$, there is an $m$ such that for all finite
  traceable graphs $G$ with $|G| \geq m$, either $G$ contains a copy
  of $K_{2,2}$ or $G$ contains a chordless $n$-path.
\end{thm}

\begin{proof}
  Let $n' = \max \{ n+1, 8 \}$ and let $m$ be the finite Ramsey number
  for $4$-tuples with $(n-1)^2+1$ many colors and a homogeneous set of
  size $n'$.  We claim that this $m$ satisfies the theorem.

  Let $G = (V,E)$ be a finite traceable graph with $|G| \geq m$.
  Without loss of generality, we assume that $V$ is an initial segment
  of $\mathbb{N}$ and that $E(i,i+1)$ holds for all $i < |G|$.  If $G$
  contains a chordless $n$-path, then we are done.  Hence assume that
  $G$ does not contain such a chordless $n$-path.

  For each $x < y$ in $G$, fix a chordless path $x = a_0(x,y) <
  a_1(x,y) < \cdots < a_{N(x,y)}(x,y) = y$.  (By assumption, there is
  an increasing path from $x$ to $y$ and hence there is a minimal
  length increasing path from $x$ to $y$ which is necessarily
  chordless.)  Note that $1 \leq N(x,y) \leq n-2$.  We use these fixed
  chordless paths to define a coloring of $[G]^4$ using $(n-1)^2+1$
  many colors.  For each $0 \leq i, j \leq n-2$ let
  \[
  K_{i,j} = \{ \{ x,y,u,v \} \mid x < y < u < v \, \wedge \, N(x,y)
  \geq i \, \wedge \, N(u,v) \geq j \, \wedge \, E(a_i(x,y),a_j(u,v))
  \}.
  \]
  That is, an increasing 4-tuple $\langle x,y,u,v \rangle$ is assigned
  color $K_{i,j}$ if the $i$-th vertex in the fixed chordless path
  from $x$ to $y$ is connected to the $j$-th vertex in the fixed
  chordless path from $u$ to $v$.  Let
  \[
  K = \{ \{ x,y,u,v \} \mid x < y < u < v \, \wedge \, \forall i \neq
  j \leq n-1 \, ( \{ x,y,u,v \} \not \in K_{i,j}) \}
  \]
  be the color of any element of $[G]^4$ not colored by any of the
  $K_{i,j}$ colors.  (Note that some $4$-tuples may be assigned more
  than one color of the form $K_{i,j}$.  This does not affect the
  statement of the finite Ramsey theorem.)

  By the finite Ramsey theorem, $G$ must have a homogeneous set of
  size $n'$ for one of these colors.  First consider the case when
  this homogeneous set is for a color $K_{i,j}$.  Since $n' \geq 8$,
  we have elements $x_1 < x_2 < \cdots < x_8$ in our homogeneous set.
  By the definition of $K_{i,j}$, we have at least the following edges
  in $G$: \bigskip

  $
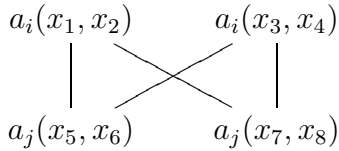
\xymatrix{ \a{d} a_i(x_1,x_2) \a{dr} & \a{d} a_i(x_3,x_4) \a{dl}  \\
    a_j(x_5,x_6) & a_j(x_7,x_8) }$

  \bigskip
  \noindent
  Thus, the vertices $a_i(x_1,x_2)$, $a_i(x_3,x_4)$, $a_j(x_5,x_6)$
  and $a_j(x_7,x_8)$ form a copy of $K_{2,2}$.

  Second, consider the case when this homogeneous set is for the color
  $K$.  In this case we derive a contradiction by showing that $G$ has
  a chordless $n$-path.  Since $n' \geq n+1$, we have elements $x_0 <
  x_1 < \cdots < x_{n-1} < x_n$ in our homogeneous set.  Thus, we have
  a path
  \begin{gather*}
    x_0 = a_0(x_0,x_1) < a_1(x_0,x_1) < \cdots <
    a_{N(x_0,x_1)}(x_0,x_1) = x_1 = a_0(x_1,x_2) <
    \cdots  \\
    \cdots < a_{N(x_1,x_2)}(x_1,x_2) = x_2 = a_0(x_2,x_3) < \cdots <
    a_{N(x_{n-1},x_n)}(x_{n-1},x_n) = x_n.
  \end{gather*}
  Define a sequence of vertices $y_0 < y_1 < \ldots$ from the vertices
  of this path as follows: take $y_0=x_0$ and for $i>0$ take $y_{i+1}$
  to be the greatest vertex $w$ on the path such that $E(y_i,w)$
  holds. Continue until either $y_{n-1}$ has been defined or until the
  vertices of the path have been exhausted. Since $\{x_0, \ldots,
  x_n\}$ is homogeneous for $K$, it follows that for each $x_j \leq v
  \leq x_{j+1}$ on the path ($0 \leq j \leq n-2$), the greatest vertex
  $w$ on the path such that $E(v,w)$ holds satisfies $w \leq x_{j+2}$,
  and so for all $i$ we have $y_i \leq x_{i+1}$. This shows that the
  process of defining the $y_i$'s terminates with the definition of
  $y_{n-1}$. By construction, $\{y_0, \ldots, y_{n-1}\}$ is the vertex
  set of a chordless path.
\end{proof}

It would be of interest to know the minimum $m=m(n)$ such that all
finite traceable graphs $G$ with $|G|\geq m$ either contain a copy of
$K_{2,2}$ or contain a chordless $n$-path. Our proof shows that there
is a constant $c>0$ such that for all $n \geq 2$ we have
$$
m(n) \leq t_{c\lceil \log n \rceil}(2)
$$
where the tower function $t_k(x)$ is defined recursively by $t_1(x) =
x$ and $t_k(x) = 2^{t_{k-1}(x)}$ for $k > 1$. (This is an easy
calculation based on known bounds for finite Ramsey numbers.)
Presumably this is far from the truth, but any substatial improvement
would require a new aproach to the proof of Theorem \ref{thm:kramsey}.

Before proving Theorem \ref{thm:GRS3} in $\RCA$, we give numerous
definitions from lattice theory within $\RCA$.  A \textit{lattice} is
a quadruple $(L, \leq_L, \wedge_L, \vee_L)$ such that $L \subseteq
\mathbb{N}$, $\leq_L$ is a binary relation on $L$ satisfying the
axioms for a partial order, and $\wedge_L$ and $\vee_L$ are functions
from $L \times L$ into $L$ such that for all $x,y \in L$, $x \wedge_L
y$ is the greatest lower bound of $x$ and $y$, and $x \vee_L y$ is the
least upper bound of $x$ and $y$.  (Typically we will drop the
subscripts on $\leq$, $\wedge$ and $\vee$.)  We denote the least
element of $L$ (if it exists) by $0_L$ and we denote the greatest
element of $L$ (if it exists) by $1_L$.

A \textit{lattice of length 3} is a lattice with a least element and a
greatest element such that every element $x \not \in \{ 0_L, 1_L \}$ is either
an atom (i.e.~there are no elements $y$ such that $0_L < y < x$) or a
coatom (i.e.~there are no elements $y$ such that $x < y < 1_L$).

\begin{lem}[$\RCA$]
  \label{lem:length3}
  Let $L$ be a lattice of length 3.  There do not exist atoms $x \neq
  y$ and coatoms $u \neq v$ such that $x <_L u$, $x <_L v$, $y <_L u$ and $y <_L v$.
\end{lem}

\begin{proof}
  Suppose for a contradiction there are such elements.  Since atoms
  are incomparable $x <_L x \vee y$ and since coatoms are incomparable
  $x \vee y <_L u$.  Therefore, $0_L <_L x <_L x \vee y <_L u <_L
  1_L$, contradicting the definition of length 3.
\end{proof}

Let $L$ be a lattice.  For each finite subset $\{ g_0, g_1, \ldots,
g_k \}$ of elements of $L$, we define an increasing sequence of finite
subsets $F_1 \subseteq F_2 \subseteq F_3 \subseteq \cdots$ of $L$ by
\begin{gather*}
  F_1 = \{ g_0, g_1, \ldots, g_k \} \\
  F_{n+1} = \{ x \wedge y \mid x,y \in F_n \} \cup \{ x \vee y \mid
  x,y \in F_n \}.
\end{gather*}
(More formally, we define a sequence of finite set codes for these
sets. Although we can form this sequence of finite sets, we cannot in
general form their union in $\RCA$ as that uses $\Sigma^0_1$
comprehension.)  We say $L$ \textit{is finitely generated} if there
exists a finite set $\{ g_0, \ldots, g_k \}$ such that $\forall x \in
L \, \exists n \, (x \in F_n)$.

A lattice $L$ \textit{contains arbitrarily long finite fences} if for
every odd $n$, there is a sequence of elements $x_0, x_1, \ldots, x_n$
of $L$ such that the Hasse diagram for these elements looks like

\bigskip

$
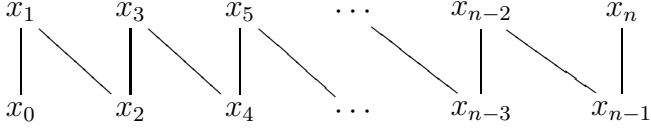
\xymatrix{ \a{d} x_1 \a{dr} & \a{d} x_3 \a{dr} & \a{d} x_5 \a{dr} &
  \cdots & \a{d} x_{n-2} \a{dr} &
  \a{d} x_{n}  \\
  x_0 & x_2 & x_4 & \cdots & \a{ul} x_{n-3} & x_{n-1} }$

\bigskip
\noindent
That is, $x_0 <_L x_1$, for each even $i$ with $0 < i < n$, $x_i <_L
x_{i-1}$ and $x_i <_L x_{i+1}$, and no other comparability relations
hold between these elements.

We can now formalize the proof of Theorem \ref{thm:GRS3} (restated
below) in $\RCA$.  The classical part of this proof is a
straightforward formalization of the proof given in \cite{GRS:82} with
an application of Theorem \ref{thm:kramsey} in place of an application
of Theorem \ref{thm:GRS1}.

\begin{thm}[$\RCA$]
  Every finitely generated infinite lattice of length 3 contains
  arbitrarily long finite fences.
\end{thm}

\begin{proof}
  Because this theorem is a $\Pi^1_1$ statement and $\RCA$ is
  conservative over $\WKL$ for $\Pi^1_1$ statements, it suffices to
  give a proof in $\WKL$.  Therefore, we work in $\WKL$.

  Let $L$ be an infinite lattice of length 3 which is finite generated
  by $\{ g_0, \ldots, g_k \}$.  Define the finite subsets $F_0
  \subseteq F_1 \subseteq \cdots$ as above.  We say that an element $x
  \in L$ has rank $0$ if $x \in F_0$ (i.e.~$x$ is a generator of $L$).
  We say $x$ has rank $n+1$ if $x \in F_{n+1} \setminus F_n$.  Note
  that every element has a rank and there are only finitely many
  elements of each rank.  Therefore, since $L$ is infinite, for every
  $n \in \mathbb{N}$, there is an element of rank $n$.  Furthermore,
  there is a function $r(x)$ giving the rank of each element and there
  is a function $m(n)$ such that
  \[
  \forall x \in L \, \forall n \in \mathbb{N} \, (r(x) = n \rightarrow
  x \leq m(n)).
  \]

  Form a tree $T \subseteq (L \setminus \{ 0_L,1_L \})^{< \mathbb{N}}$
  as follows.  The sequence $\langle x_0, x_1, \ldots, x_n \rangle \in
  T$ if and only if for every $i \leq n$, $r(x_i) = i$, and for every
  $0 < i \leq n$, there is an $a \in L$ with $r(a) < i$ such that $x_i
  = x_{i-1} \vee a$ or $x_i = x_{i-1} \wedge a$.  $T$ has the
  following properties.
  \begin{enumerate}
  \item[(P1)] If $\langle x_0, \ldots, x_n \rangle \in T$, then the
    $x_i$ are distinct and for all $i < n$, $x_i$ is comparable with
    $x_{i+1}$.  This property follows since $r(x_i) = i$, $r(x_{i+1})
    = i+1$ and $x_{i+1} = x_i \wedge a$ or $x_{i+1} = x_i \vee a$ for
    some $a \in L$.  (Because $L$ has length 3, the $x_i$ elements
    alternate between atoms and coatoms.)

  \item[(P2)] For every $x \in L \setminus \{ 0_L, 1_L \}$, there is a
    $\sigma \in T$ such that $\sigma*x \in T$.  This property follows
    by induction on the rank of $x$.  If the $r(x) = 0$, then $\sigma
    = \emptyset$.  If $r(x) = n+1$, then $x \in F_{n+1} \setminus
    F_n$.  Without loss of generality suppose $x = y \wedge z$ where
    $y,z \in F_n$.  Then $r(y), r(z) \leq n$ and either $r(y) = n$ or
    $r(z) = n$.  Suppose $r(y) = n$.  By the induction hypothesis,
    there is a $\tau \in T$ such that $\tau*y \in T$.  Let $\sigma =
    \tau*y$ and it follows from the definition of $T$ that $\sigma*x
    \in T$.

  \item[(P3)] $T$ is infinite.  This property follows from (P2) and
    the fact that $L$ has elements of rank $n$ for each $n \in
    \mathbb{N}$.

  \item[(P4)] The branching in $T$ is bounded by the function $m(n)$
    in the sense that if $\sigma \in T$ then for all $i < |\sigma|$,
    $\sigma(i) \leq m(i)$.  This property follows from the definition
    of $T$.
  \end{enumerate}

  Since $T$ is an infinite tree with bounded branching, $\WKL$ proves
  that $T$ has an infinite path $f: \mathbb{N} \rightarrow L \setminus
  \{ 0_L, 1_L \}$.  (See Lemma IV.1.4 in Simpson \cite{sim:book}.)  By
  (P1), $f$ is 1-to-1.  Furthermore the range of $f$ exists since $x
  \in \text{range}(f)$ if and only if $f(n) = x$ where $n = r(x)$.  If
  $f(0)$ is an atom, then the Hasse diagram of the range of $f$
  contains at least the following comparability relations \bigskip

  $
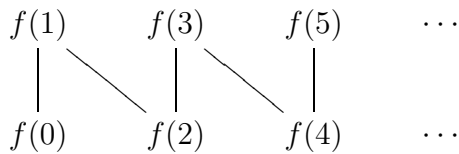
\xymatrix{ \a{d} f(1) \a{dr} & \a{d} f(3) \a{dr} & \a{d} f(5) & \cdots & \\
    f(0) & f(2) & f(4) & \cdots & }$

  \bigskip
  \noindent
  and may contain additional comparability relations.  If $f(0)$ is a
  coatom, then we obtain the dual of this picture.  To avoid breaking
  into simple dual cases, we will assume $f(0)$ is an atom for the
  remainder of the proof.

  Define a graph $G = (V,E)$ with $V = \text{range}(f)$ and
  $E(f(n),f(m))$ holds if and only if $f(n)$ and $f(m)$ are comparable
  in $L$.  $G$ looks like the Hasse diagram above with possibly
  additional edges (since the lattice elements in this diagram could
  have additional comparability relations).  However, each $f(2n)$ is
  an atom and each $f(2n+1)$ is a coatom.  Therefore, by Lemma
  \ref{lem:length3}, $G$ does not contain a copy of $K_{2,2}$.  By
  (P1), $f$ is a tracing function for $G$, so $G$ is an infinite
  traceable graph that does not contain a copy of $K_{2,2}$.
  Therefore, by Theorem \ref{thm:kramsey}, $G$ contains arbitrarily
  long finite chordless paths.  Since finite chordless paths in $G$
  are finite fences when viewed in $L$, $L$ contains arbitrarily long
  finite fences.

\end{proof}

\end{document}